\documentclass[a4paper,12pt]{article}

\usepackage{amsmath}
\usepackage{amsthm}
\usepackage{amssymb}
\usepackage{color}

\def\R{\mathbb R}
\def\Z{\mathbb Z}

\def\tu{\tilde{u}}

\def\tb{\tilde{b}}
\def\tX{\tilde{X}}
\def\tsigma{\tilde{\sigma}}

\def\M{\mathcal{M}}

\def\epsilon{\varepsilon}
\def\ds{\displaystyle}

\newcommand{\be}{\begin{equation}}
\newcommand{\ee}{\end{equation}}
\newcommand{\baa}{\begin{array}}
\newcommand{\eaa}{\end{array}}
\newcommand{\ba}{\begin{eqnarray}}
\newcommand{\ea}{\end{eqnarray}}
\newcommand{\bi}{\begin{itemize}}
\newcommand{\ei}{\end{itemize}}

\newtheorem{lemma}{Lemma}[section]
\newtheorem{theorem}[lemma]{Theorem}
\newtheorem{corollary}[lemma]{Corollary}

\newtheorem{remark}[lemma]{Remark}

\begin{document}

\title{Simultaneous determination of the drift and diffusion coefficients in stochastic differential equations}

\author{M. Cristofol \thanks{Institut de Math\'{e}matiques de Marseille, CNRS, UMR 7373, \'Ecole Centrale, Aix-Marseille Universit\'e, 13453 Marseille, France,
  e-mail: \texttt{\ michel.cristofol@univ-amu.fr} }
\and
L. Roques \thanks{DBioSP, INRA, 84914, Avignon, France, e-mail: \texttt{\ lionel.roques@inra.fr}}}

\date{}

\maketitle

%
%
%

%
%
%
%
%

\begin{abstract}
In this work, we consider a one-dimensional It\^o diffusion process $X_t$ with possibly nonlinear drift and diffusion coefficients. We show that, when the diffusion coefficient is known, the drift coefficient is uniquely determined by an observation of the expectation of the process during a small time interval, and starting from values $X_0$ in a given subset of $\R$. With the same type of observation, and given the drift coefficient, we also show that the diffusion coefficient is uniquely determined. When both coefficients are unknown, we show that they are simultaneously uniquely determined by the observation of the expectation and variance of the process, during a small time interval, and starting again from values $X_0$ in a given subset of $\R$. To derive these results, we apply the Feynman-Kac theorem which leads to a linear parabolic equation with unknown coefficients in front of the first and second order terms. We then solve the corresponding inverse problem with PDE technics which are mainly based on the strong parabolic maximum principle.
\end{abstract}


\section{Introduction}

We consider one-dimensional It\^o diffusion processes $X_t \in \R$ satisfying  stochastic differential equations of the form:
\be
d X_t= b(X_t) \,  dt + \sigma(X_t) dW_t, \ t\in[0,T]; \ X_0=x,
\label{eq:itodiffusion}
\ee
where $T>0$, $W_t$ is the one-dimensional Wiener process and $b: \, \R \to \R,$ $\sigma: \, \R \to \R, \sigma > 0, $ are Lipschitz-continuous functions. Under these assumptions, the solution of the equation \eqref{eq:itodiffusion} is unique in the sense of theorem 5.2.1 in \cite{Oks03}. The term $b(X_t)\, dt $ can be interpreted as the deterministic part of the equation, while $\sigma(X_t) dW_t$ is the stochastic part of the equation. In the sequel, the functions $b$ and $\sigma$ are called the \emph{drift term} and \emph{diffusion term}, respectively.

These equations arise in several domains of applications, such as biology, physics and financial mathematics. We detail below some classical forms of the functions $b$ and $\sigma$: (1) in Ornstein-Uhlenbeck processes,  $b(X_t)=\theta \, (\mu-X_t)$ and $\sigma(X_t)=\sigma=\hbox{cte}$ with $\mu \in \R$ and $\theta,\sigma>0$. Ornstein-Uhlenbeck processes describe a noisy relaxation with equilibrium $\mu$. They find applications in physics \cite{VanHen92}, financial mathematics \cite{NicVen03}
 and biology \cite{SmoFoc10}; (2) in the two types Wright-Fisher gene frequency diffusion model with selection and genetic drift effects, $b(X_t)=X_t\left[m_1-(m_1\,X_t +m_2 \, (1-X_t) )\right]$ and $\sigma(X_t)=\sqrt{\frac{1}{N_e} \, X_t(1-X_t)}$, for some constants $m_1, \, m_2, \, N_e$; this is one of the most standard model in population genetics \cite{Ewe12}; (3) in Geometric Brownian motion, $b(X_t)=\alpha X_t$ and $\sigma(X_t)=\beta X_t$.  This equation is used in finance, with non-constant coefficients $\alpha$, $\beta$,  to model stock prices in the Black-Scholes model. The term  $\alpha$ is interpreted in this case as the percentage drift and $\beta$ the percentage volatility \cite{Hul06}.
  The determination of the volatility is an important question in finance, and is generally addressed numerically based on observations of the prices of financial options \cite{LisYou01,DenYu08}; see also \cite{BouIsa99} for a uniqueness result based on the same type of observations.

The aim of our study is to determine the drift term $b$ and the diffusion term $\sigma$ for general equations of the form~\eqref{eq:itodiffusion}, based on observations of the stochastic process $X_t$. Equivalently, this means showing the uniqueness of the coefficients $b$ and $\sigma$ which lead to a solution that matches with the given observation. The main type of observation that we consider is the expectation $E^{x}[f(X_t)]=E[f(X_t)|X_0=x]$, of some function of the stochastic process $X_t,$ for instance a momentum if $f(s)=s^k$ for some $k\ge 0$. The observation is carried out during a small time interval and for initial conditions $X_0$ in a small subset of $\R$. In that respect we use parabolic partial differential equation (PDE) technics inspired from the theory of inverse problems.

The It\^o diffusion processes are related to PDEs by the Feynman-Kac theorem (see e.g. theorem 8.1.1 in \cite{Oks03}). Consider a function
 \be \label{hyp:f_FK} f\in C^2(\R) \hbox{ such that }|f(x)|\le C \, e^{\delta x^2},\ee for $\delta>0$ small enough and some $C>0$. Define
\be
u(t,x)=E^x\left[f(X_t)\right]=E\left[f(X_t)|X_0=x\right],
\ee
where $X_t$ is the solution of \eqref{eq:itodiffusion} with $X_0=x$. The Feynman-Kac theorem implies that $u$ is the unique solution in $C^2_1(\R_+ \times \R)$ of:
\be \label{eq:para_gale}
\partial_t u =\frac{1}{2} \sigma^2(x) \partial_{xx} u+ b(x) \partial_x u,  \ t\ge 0; \ u(0,x)=f(x).
\ee
For parabolic equations of the form \eqref{eq:para_gale}, several inverse problems have already been investigated. In all cases, the main question is to show the uniqueness of some coefficients in the equation, based on exact observations of the solution $u(t,x)$, for $(t,x)$ in a given observation region $\mathcal{O}\subset [0,+\infty)\times \R.$ Furthermore, one of the most challenging  goal is to obtain such uniqueness results using the smallest possible observation region.

Most uniqueness results in inverse problems for parabolic PDEs have been obtained using the method of Carleman estimates \cite{BukKli81} on bounded domains. This method requires,
among other measurements, the knowledge of the solution $u(\tau,x)$ at some time $\tau>0$ and \emph{for all} $x$ in the domain \cite{YamZou01,BelYam06,CriRoq08,Yam09,ImmYam98,CriRoq13}. Other approaches are based on a semi-group formulation of the solutions, but use the same type of observations of the solution on the whole domain, at a given time \cite{ChoYam08}.
More recent approaches \cite{RoqCri10,CriGarHamRoq11,RoqCri12,RoqCheCriSouGhi14} lead to uniqueness results for one or several coefficients, under the assumption that $u$ and its first spatial
derivative are known at a \textit{single point} $x_0$ of a bounded domain, and for all $t$ in a small interval $(0,\varepsilon),$ and
that the initial data $u(0,x)$  is known over the entire domain. On the other hand, the case of unbounded domains  is less addressed (see \cite{CriKad14}).\\
Here, contrarily to most existing approaches, (i) the domain is unbounded; (ii) we determine simultaneously two coefficients in front of a second and a first order term in the PDE; (iii) our results are interpreted in terms of nonlinear stochastic diffusion processes. As in the above-mentioned studies \cite{RoqCri10,CriGarHamRoq11,RoqCri12,RoqCheCriSouGhi14}, we assume that the observation set reduces to a neighborhood of single point $x_0$, during a small time interval.

Our manuscript is organized as follows. In Section~\ref{sec:assum_res} we detail our assumptions on the unknown coefficients and on the observations and we state our main results. In Section~\ref{sec:proofs} we prove the uniqueness results stated in Theorem \ref{th:1coeff} , Theorem \ref{th:scoeff} and Theorem \ref{th:2coeffs}.

\section{Assumptions and main results \label{sec:assum_res}}

\paragraph{Observations.} We consider two main types of observations. Let $\varepsilon \in (0,T)$ and $\omega$ an open and nonempty subset in $\R$.  The observation sets are either of the form
\be \label{obs1}
\mathcal{O}_f[X_t]=\{E^{x}[f(X_t)], \hbox{ for }t\in (0,\varepsilon) \hbox{ and }x\in \omega\},
\ee
for some function $f$ satisfying the assumptions \eqref{hyp:f_FK} of Feynman-Kac theorem, or of the form
\be
\label{obs2}
\mathcal{O}^k[X_t]=\{E^{x}[(X_t)^k], \hbox{ for }t\in (0,\varepsilon) \hbox{ and }x\in \omega\},
\ee
for  $k=1,2$. In both cases,  $\varepsilon>0$ and $\omega$ can be chosen as small as we want. For the sake of simplicity, and with a slight abuse of notation, for two processes $X$ and $\tX$, we say that $\mathcal{O}_f[X_t]=\mathcal{O}_f[\tX_t]$ (resp. $\mathcal{O}^k[X_t]=\mathcal{O}^k[\tX_t]$) if and only if $E^{x}[f(X_t)]=E^{x}[f(\tX_t)]$ (resp. $E^{x}[(X_t)^k]=E^{x}[(\tX_t)^k]$ for $k=1,2$), for all $t\in (0,\varepsilon)$ and $x\in \omega$.

\

\paragraph{Unknown functions.} We assume that the unknown functions belong to the function space:
\be
\label{hyp:piecewise}
\M:=\{\psi \hbox{ is Lipschitz-continuous and piecewise analytic in }\R\}.
\ee
A continuous function $\psi$ is called piecewise analytic
if there exist $n \ge 1$ and an increasing sequence $(\kappa_j)_{j\in \Z}$ such that $\lim\limits_{j\to -\infty}\kappa_j=-\infty$, $\lim\limits_{j\to +\infty}\kappa_j=+\infty$, $\kappa_{j+1}-\kappa_{j}>\delta$ for some $\delta>0$, and  \\
$\ds{\psi(x)=\sum_{j\in\Z} \chi_{[\kappa_j,\kappa_{j+1})}(x)\varphi_j(x)},$
for all $x\in\R;$ \\
  here $\varphi_j$ are some analytic functions defined on the intervals
$[\kappa_j,\kappa_{j+1}]$, and
$\chi_{[\kappa_j,\kappa_{j+1})}$ are the characteristic functions
of the intervals $[\kappa_j,\kappa_{j+1})$ for $j\in\Z$.

In practice, the assumption $\psi \in \M$ is not very restrictive. For instance, the set of piecewise linear functions in $\R$ is a subset of $\M.$

\

\paragraph{Main results.} Our first result states that, whenever $\sigma$ is known, the coefficient $b$ in \eqref{eq:itodiffusion} is uniquely determined by an observation of the type $\mathcal{O}_f.$
\begin{theorem}\label{th:1coeff} Let $b$ and $\tb$ in $\M$, $\sigma$ a strictly positive Lipschitz-continuous function, $X_t$ the solution of \eqref{eq:itodiffusion}, and $\tX_t$ the solution of $d \tX_t= \tb(\tX_t) \,  dt + \sigma(\tX_t) dW_t, \ t\in[0,T]; \ \tX_0=x$. Assume  that $f'\neq 0$ in $\R$ and $\mathcal{O}_f[X_t]=\mathcal{O}_f[\tX_t]$. Then, $b\equiv \tb$ in $\R$.
\end{theorem}
An important and easily interpretable observation is the expectation of the process $X_t$ during a small time interval and for all $X_0=x$ in any small set $\omega \subset \R$ by choosing $f(x) = x, \; x \in \R$.\\

Our second result allows to uniquely determine, whenever $b$ is known,  the coefficient $\sigma$ in \eqref{eq:itodiffusion}, which is a coefficient from the principal part (second order term) of the equation \eqref{eq:para_gale}.
\begin{theorem}\label{th:scoeff} Let $\sigma$ and $\tsigma >0$ in $\M$, $b$ a Lipschitz-continuous function, $X_t$ the solution of \eqref{eq:itodiffusion}, and $\tX_t$ the solution of $d \tX_t= b(\tX_t) \,  dt + \tsigma(\tX_t) dW_t, \ t\in[0,T]; \ \tX_0=x$. Assume  that $f'' \neq 0$ in $\R$ and $\mathcal{O}_f[X_t]=\mathcal{O}_f[\tX_t]$. Then, $\sigma \equiv \tsigma$ in $\R$.
\end{theorem}

Determining several coefficients of parabolic PDEs is generally far more involved than determining a single coefficient. It requires more and well-chosen observations. For instance, four coefficients of a Lotka-Volterra system of parabolic equations have been determined in \cite{RoqCri12}, based on the observation of one component of the solution, starting with three different initial conditions. See also \cite{CriRoq13,ChoYam08} for other results on simultaneous determination of several coefficients, with different methods. Here, our third result shows that, if the first momentum (expected value) and the second momentum of $X_t$ are observed during a small time interval and for $X_0=x$ in a small set $\omega \subset \R$, then  both coefficients $b$ and $\sigma$  in \eqref{eq:itodiffusion} are uniquely determined.
\begin{theorem}\label{th:2coeffs} Let $b,\tb ,\sigma,\tsigma \in \M$ with $\sigma,\tsigma>0$. Consider $X_t$ the solution of \eqref{eq:itodiffusion}, and $\tX_t$ the solution of  $d \tX_t= \tb(\tX_t) \,  dt + \tsigma(\tX_t) dW_t, \ t\in[0,T]; \ \tX_0=x$. Assume that $\mathcal{O}^k[X_t]=\mathcal{O}^k[\tX_t]$  for $k=1,2$. Then, $b\equiv \tb$ and $\sigma\equiv\tsigma$ in $\R$.
\end{theorem}
An immediate corollary of Theorem~\ref{th:2coeffs} is that $b$ and $\sigma$  are uniquely determined by the observation of the expectation $E^{x}[X_t]$ and variance $V^x[X_t]=E^{x}[X_t^2]-(E^{x}[X_t])^2$ of the process $X_t$ during a small time interval and for all $X_0=x$ in a small set $\omega \subset \R$. More precisely, define the set
\be
\mathcal{O}_v[X_t]=\{V^x[X_t], \hbox{ for }t\in (0,\varepsilon) \hbox{ and }x\in \omega\},
\ee
we have the following result.
\begin{corollary}\label{th:cor1} Let $b,\tb ,\sigma,\tsigma\in \M$. Consider $X_t$ the solution of \eqref{eq:itodiffusion}, and $\tX_t$ the solution of $d \tX_t= \tb(\tX_t) \,  dt + \tsigma(\tX_t) dW_t, \ t\in[0,T]; \ \tX_0=x,$ respectively. Assume that $\mathcal{O}^1[X_t]=\mathcal{O}^1[\tX_t]$ and $\mathcal{O}_v[X_t]=\mathcal{O}_v[\tX_t]$. Then, $b\equiv \tb$ and $\sigma\equiv \tsigma$ in $\R$.
\end{corollary}

\begin{remark}
All of our results remain true if the observations \eqref{obs1} and \eqref{obs2} are replaced by pointwise observations at a given point $x_0\in \R$ instead of observations in a subdomain $\omega$. More precisely, if  \eqref{obs1} and \eqref{obs2} are replaced by
\be \label{obs1b}
\mathcal{O}^{'}_f[X_t]=\{E^{x_0}[f(X_t)], \, \partial_x E^{x}[f(X_t)]_{|x=x_0}, \hbox{ for }t\in (0,\varepsilon)\},
\ee
and
\be
\label{obs2b}
\mathcal{O}^{',k}[X_t]=\{E^{x_0}[(X_t)^k], \,  \partial_x E^{x}[(X_t)^k]_{|x=x_0}, \hbox{ for }t\in (0,\varepsilon)\},
\ee
for  $k=1,2$, all of the results of our theorems and corollary can still be obtained, by using the Hopf's Lemma in addition to the strong parabolic maximum principle. See the footnote in the proof of Theorem~\ref{th:1coeff}.
\end{remark}

\section{Proofs \label{sec:proofs}}

\subsection{Proof of Theorems \ref{th:1coeff} and \ref{th:scoeff}}

We begin with the proof of Theorem~\ref{th:1coeff}.
Let us define, for all $t\in[0,T)$ and $x\in \R$,
\be \baa{l}
\ds u(t,x)=E^x\left[f(X_t)\right]=E\left[f(X_t)|X_0=x\right], \\ \ds \tu(t,x)=E^x[f(\tX_t)]=E[f(\tX_t)|X_0=x].\eaa
\ee
As mentioned in the Introduction, the Feynman-Kac theorem implies that $u$ and $\tu$ are respectively the unique solutions of:
\be \label{eq:para_u}
\partial_t u =\frac{1}{2} \sigma^2(x) \partial_{xx} u+ b(x) \partial_x u,  \ t\in[0,T), \ x\in \R; \ u(0,x)=f(x),
\ee
and
\be \label{eq:para_tu}
\partial_t \tu =\frac{1}{2} \sigma^2(x) \partial_{xx} \tu+ \tb(x) \partial_x \tu,  \ t\in[0,T), \ x\in \R \ \tu(0,x)=f(x).
\ee
Define $$B(x)=b(x)-\tb(x), \hbox{ and } U(t,x)=u(t,x)-\tu(t,x).$$
Then $U(t,x)$ satisfies
\be \label{eq:U}
\partial_t U =\frac{1}{2} \sigma^2(x) \partial_{xx} U+ b(x) \partial_x U +
 B(x) \, \partial_x \tu,  \ t\in[0,T), \ x\in \R,
\ee
and $U(0,x)=0$ for all $x\in \R$.

Let $x_0\in \omega$. As $B\in \M$ is piecewise analytic, we can define
$$x_1=\sup\{x>x_0 \hbox{ such that } B \hbox{ has a constant sign on }[x_0,x]\}.$$ By ``constant sign", we mean that either $B\ge 0$ on $[x_0,x]$ or $B\le 0$ on $[x_0,x]$.

Assume (by contradiction) that there exists $x_2\in (x_0,x_1)$ such that $|B(x_2)|>0$. From the definition of $x_1$, we know that $B$ has a constant sign in $(x_0,x_2)$. As $\partial_x \tu(0,x)=f'(x)\neq 0$ on the compact set $[x_0,x_2]$ and from the regularity of $\tu$ (see theorem 8.2.1 in \cite{Oks03}), there exists $\varepsilon'\in (0,\varepsilon)$ such that  $\partial_x \tu(t,x)$ has a constant sign on $[0,\varepsilon')\times [x_0,x_2]$. Finally, $B(x) \, \partial_x \tu,$ has a constant sign in $(0, \varepsilon')\times [x_0,x_2]$. Without loss of generality, we can assume that:
\be \label{eq:ineg_B1}
 B(x) \, \partial_x \tu \ge 0 \hbox{ for } (t,x) \hbox{ in } [0, \varepsilon')\times [x_0,x_2].
\ee
Computing \eqref{eq:U} at $t=0$ and $x=x_2$, and using the equality $U(0,x)=0$  we get $\partial_t U(0,x_2)=B(x_2) \, \partial_x \tu(0,x_2)\ge0$ (from \eqref{eq:ineg_B1}). Besides, from the assumption $|B(x_2)|>0$ and $f'(x_2)\neq 0$, we know that the inequality is strict: $\partial_t U(0,x_2)=B(x_2) \, f'(x_2)>0$. Thus, (even if it means reducing $\varepsilon'>0$),  \be\label{eq:ineg_B2}U(t,x_2)>0 \hbox{ for }t\in(0,\varepsilon'). \ee
Using the assumption $\mathcal{O}_f[X_t]=\mathcal{O}_f[\tX_t]$ of Theorem~\ref{th:1coeff}, and from the definition of $u$, $\tu$ and $U=u-\tu$, we have:
\be \label{eq:ineg_B3}
U(t,x)\equiv 0 \hbox{ in }[0,\varepsilon)\times \omega.
\ee
In particular, $U(t,x_0)=0$ for all $t\in (0,\varepsilon' )$. Setting $\mathcal{L}U=\frac{1}{2} \sigma^2(x) \partial_{xx} U+ b(x) \partial_x U,$ and summarizing the properties \eqref{eq:ineg_B1}-\eqref{eq:ineg_B3}, we get:
\be
\left\{\baa{l}
\ds \partial_t U - \mathcal{L} U \ge 0, \ t\in (0, \varepsilon'),  \  x \in (x_0,x_2),\vspace{0.1cm}\\
\ds U(t,x_0)=0, \ U(t,x_2)>0, \ t\in (0, \varepsilon'), \vspace{0.1cm}\\
\ds U(0,x)=0,  \  x \in (x_0,x_2).
\eaa \right.
\label{eq:U_sys}
\ee
The strong parabolic maximum principle then implies that  $U(t,x)>0$ in $(0,\varepsilon')\times(x_0,x_2).$ This contradicts \eqref{eq:ineg_B3}\footnote{If \eqref{eq:ineg_B3} was replaced by $U(t,x_0)=\partial_x U(t,x_0)= 0$  for $t\in [0,\varepsilon),$ a similar contradiction could be obtained by using the Hopf's Lemma (theorem 14 p. 190 in \cite{ProWei67}), as it implies that $\partial_x U(t,x_0)\neq 0.$}; as a consequence, $B\equiv 0$ in $(x_0,x_1)$. From the definition of $x_1$ and the piecewise analyticity of $B$, this implies that $x_1=+\infty$, thus $B\equiv 0$ in $(x_0,+\infty)$. Using the same arguments with $x_1^-=\inf\{x<x_0 \hbox{ such that } B \hbox{ has a constant sign on }[x,x_0]\}$ instead of $x_1$, we easily see that $B\equiv 0$ in $(-\infty,x_0)$, and consequently, $B\equiv 0$ in $\R$. This concludes the proof of Theorem~\ref{th:1coeff}. $\Box$

The proof of Theorem~\ref{th:scoeff} is very similar to that of Theorem~\ref{th:1coeff}.

\subsection{Proof of Theorem \ref{th:2coeffs}}
In this case, the proof is more involved. Indeed, we reconstruct simultaneously  two coefficients from the  principal part and the first order term in equation \eqref{eq:para_gale} and this implies to repeat the observations and to consider adapted weight functions in the form \eqref{obs2}. As in the proof of Theorem \ref{th:1coeff}, we define, for all $t\in[0,T)$ and $x\in \R$, and for $f(s)=s^k$,
\be \baa{l}
\ds u(t,x)=E^x[f(X_t)]=E[f(X_t)|X_0=x], \\ \ds \tu(t,x)=E^x[f(\tX_t)]=E[f(\tX_t)|X_0=x],\eaa
\ee
and $u$ and $\tu$ are respectively the unique solutions of:
\be \label{eq:para_u2}
\partial_t u =\frac{1}{2} \sigma^2(x) \partial_{xx} u+ b(x) \partial_x u,  \ t\in[0,T), \ x\in \R; \ u(0,x)=f(x),
\ee
and
\be \label{eq:para_tu2}
\partial_t \tu =\frac{1}{2} \tsigma^2(x) \partial_{xx} \tu+ \tb(x) \partial_x \tu,  \ t\in[0,T), \ x\in \R \ \tu(0,x)=f(x).
\ee
Define $$B(x)=b(x)-\tb(x), \ \Sigma(x)=\frac{1}{2}(\sigma^2(x)-\tsigma^2(x)), \hbox{ and } U(t,x)=u(t,x)-\tu(t,x).$$
Then $U(t,x)$ satisfies
\be \label{eq:U2}
\partial_t U - \frac{1}{2} \sigma^2(x) \partial_{xx} U - b(x) \partial_x U =
 B(x) \, \partial_x \tu + \Sigma(x) \, \partial_{xx} \tu,  \ t\in[0,T), \ x\in \R,
\ee
and $U(0,x)=0$ for all $x\in \R$.

Let $x_0\in \omega$. We define:
\be
\baa{l}
x^*_B=\sup\{x>x_0 \hbox{ such that } B\equiv 0 \hbox{ on }[x_0,x]\}, \\
x^*_\Sigma=\sup\{x>x_0 \hbox{ such that } \Sigma \equiv 0 \hbox{ on }[x_0,x]\}.
\eaa
\ee
Then, four cases may occur.

\

\noindent \textit{Case 1: we assume that $x^*_B<x^*_\Sigma$.} Using the piecewise analyticity of $B$, and from the definition of $x^*_B$, we obtain the existence of some $x_2\in (x^*_B,x^*_\Sigma)$ such that $B(x)\neq 0$ for all $x\in (x^*_B,x_2]$, i.e., $B$ has a constant strict sign in $(x^*_B,x_2]$. Moreover, $\Sigma(x)=0$ for all  $x\in (x^*_B,x_2]$, thus $U$ satisfies:
\be \label{eq:U2_case1}
\partial_t U - \mathcal{L} U = B(x) \, \partial_x \tu,  \ t\in[0,T), \ x\in (x_0,x_2),
\ee
where $ \mathcal{L} U := \frac{1}{2} \sigma^2(x) \partial_{xx} U+ b(x) \partial_x U$. Take $k=1$ in the definition of $f(s)=s^k$. We have $\partial_x \tu(0,x)=f'(x)=1$ for all $x\in \R$, which implies that there exists $\varepsilon'\in (0,\varepsilon)$ such that  $\partial_x \tu(t,x)$ is positive on $[0,\varepsilon')\times [x_0,x_2]$. Finally, the term $B(x) \, \partial_x \tu$ in the right hand side of \eqref{eq:U2_case1} has a constant sign in $(0, \varepsilon')\times [x_0,x_2]$. Without loss of generality, we can assume that:
\be \label{eq:ineg_B1_case1}
 B(x) \, \partial_x \tu \ge 0 \hbox{ for } (t,x) \hbox{ in } [0, \varepsilon')\times [x_0,x_2].
\ee
We then observe that $\partial_t U(0,x_2)=B(x_2) \ge 0$ and, from the definition of $x_2$,  the inequality is strict: $\partial_t U(0,x_2)>0$. Thus, (even if it means reducing $\varepsilon'>0$),  \be\label{eq:ineg_B2_case1}U(t,x_2)>0 \hbox{ for }t\in(0,\varepsilon'). \ee
Finally, $U$ satisfies
\be
\left\{\baa{l}
\ds \partial_t U - \mathcal{L} U \ge 0, \ t\in (0, \varepsilon'),  \  x \in (x_0,x_2),\vspace{0.1cm}\\
\ds U(t,x_0)=0, \ U(t,x_2)>0, \ t\in (0, \varepsilon'), \vspace{0.1cm}\\
\ds U(0,x)=0,  \  x \in (x_0,x_2).
\eaa \right.
\label{eq:U_sys2}
\ee
From strong parabolic maximum principle  $U(t,x)>0$ in $(0,\varepsilon')\times(x_0,x_2).$ This contradicts the assumption $\mathcal{O}^1[X_t]=\mathcal{O}^1[\tX_t]$ of Theorem~\ref{th:2coeffs}. Thus, Case 1 is ruled out.

\

\noindent \textit{Case 2: we assume that $x^*_B>x^*_\Sigma$.} With the same type of arguments as in Case 1, we obtain the existence of some $x_2\in (x^*_\Sigma,x^*_B)$ such that $\Sigma(x)\neq 0$ for all $x\in (x^*_\Sigma,x_2]$, i.e., $\Sigma$ has a constant strict sign in $(x^*_\Sigma,x_2]$. Moreover, $B(x)=0$ for all  $x\in (x^*_\Sigma,x_2]$, thus $U$ satisfies:
\be \label{eq:U2_case2}
\partial_t U - \mathcal{L} U =  \Sigma(x) \, \partial_{xx} \tu,  \ t\in[0,T), \ x\in (x_0,x_2).
\ee
Take $k=2$ in the definition of $f(s)=s^k$. We have $\partial_{xx} \tu(0,x)=f''(x)=2$ for all $x\in \R$. Thus, with the same arguments as in Case 1, we get:
\be \label{eq:ineg_B1_case2}
 \Sigma(x) \, \partial_{xx} \tu \ge 0 \hbox{ for } (t,x) \hbox{ in } [0, \varepsilon')\times [x_0,x_2],
\ee
and $\partial_t U(0,x_2)>0$. Thus, $U$ again satisfies \eqref{eq:U_sys2}, and the strong parabolic maximum principle implies $U(t,x)>0$ in $(0,\varepsilon')\times(x_0,x_2),$ leading to a contradiction with the assumption $\mathcal{O}^2[X_t]=\mathcal{O}^2[\tX_t]$ of Theorem~\ref{th:2coeffs}. Thus, Case 2 is ruled out.

\

\noindent \textit{Case 3: we assume that $x^*_B=x^*_\Sigma<+\infty$.} Let us set
$$G(t,x)=B(x) \, \partial_x \tu + \Sigma(x) \, \partial_{xx} \tu,$$corresponding to the right-hand side in \eqref{eq:U2}. Then, set $$l^*=\lim\limits_{x\to x^*_B, x> x^*_B}\frac{\Sigma(x)}{B(x)}.$$From the analyticity of $\Sigma$ and $B$ in a right neighborhood of $x^*_B$, $l^*$ is well-defined and only two situations may occur: either $|l^*|<+\infty$ or $|l^*|=+\infty$.

Assume first that $|l^*|<+\infty$. Take $k=1$ in the definition of $f(s)=s^k$. Thus,
\be\label{eq:dxdxxu1} \partial_x \tu(0,x)=1 \hbox{ and }\partial_{xx} \tu(0,x)=0.\ee  Let $x_2> x^*_B$ such that $B(x)\neq 0$ in $(x^*_B,x_2]$ and $\Sigma(x)/B(x)$ remains bounded in $(x^*_B,x_2]$. Without loss of generality, we can assume that $B>0$ in $(x^*_B,x_2]$. Using \eqref{eq:dxdxxu1}, and since $|l^*|<+\infty$, we obtain the existence of $\varepsilon'\in (0,\varepsilon)$ such that
$$\frac{G(t,x)}{B(x)}=\partial_x \tu + \frac{\Sigma(x)}{B(x)} \, \partial_{xx} \tu >0 \hbox{ for } (t,x) \hbox{ in } (0, \varepsilon')\times (x_0,x_2),$$and $G(t,x)$ satisfies the same inequality. Thus, again, $U$ satisfies \eqref{eq:U_sys2}, and the strong parabolic maximum principle implies that $U(t,x)>0$ in $(0,\varepsilon')\times(x_0,x_2)$ and a contradiction with the assumption $\mathcal{O}^1[X_t]=\mathcal{O}^1[\tX_t]$ of Theorem~\ref{th:2coeffs}. The assumption $|l^*|<+\infty$ is then ruled out.

Assume now that $|l^*|=+\infty$.  Take $k=2$ in the definition of $f(s)=s^k$. This time,
\be\label{eq:dxdxxu2} \partial_x \tu(0,x)=2\, x \hbox{ and }\partial_{xx} \tu(0,x)=2.\ee  Let $x_2> x^*_\Sigma$ such that $$\Sigma(x)\neq 0 \hbox{ and }|2\, x \, (B(x)/\Sigma(x))|<1  \hbox{ in }(x^*_\Sigma,x_2].$$ Without loss of generality, we assume that $\Sigma>0$ in $(x^*_\Sigma,x_2]$. Using \eqref{eq:dxdxxu2}, and since $|l^*|=+\infty$, we can define $\varepsilon'\in (0,\varepsilon)$ such that
$$\frac{G(t,x)}{\Sigma(x)}= \frac{B(x)}{\Sigma(x)} \partial_x \tu + \, \partial_{xx} \tu >0 \hbox{ for } (t,x) \hbox{ in } (0, \varepsilon')\times (x_0,x_2).$$Again, using the strong parabolic maximum principle, we get a contradiction with the assumption $\mathcal{O}^2[X_t]=\mathcal{O}^2[\tX_t]$ of Theorem~\ref{th:2coeffs}. Case 3 is then ruled out.

Finally, as Cases 1, 2, 3 are ruled out, we necessarily have $x^*_B=x^*_\Sigma=+\infty$, which show that $B\equiv\Sigma\equiv 0$ in $(x_0,+\infty)$. Using the same arguments with
$(x^*_B)^-=\inf\{x<x_0 \hbox{ such that } B\equiv 0 \hbox{ on }[x,x_0]\}$ and  $(x^*_\Sigma)^-=\inf\{x<x_0 \hbox{ such that } \Sigma \equiv 0 \hbox{ on }[x,x_0]\},$
 instead of $x^*_B$ and $x^*_\Sigma$, we also check that $B\equiv\Sigma\equiv 0$ in $(-\infty,x_0)$ and consequently $B\equiv\Sigma\equiv 0$  in $\R$ which concludes the proof of Theorem~\ref{th:2coeffs}. $\Box$

\section*{Acknowledgements}
The research leading to these results has received funding from the ANR within the project NONLOCAL ANR-14-CE25-0013.

\end{document}